\documentclass[11pt]{amsart}
\usepackage{amsmath}
\usepackage{amsxtra}
\usepackage{amscd}
\usepackage{amsthm}
\usepackage{amsfonts}
\usepackage{amssymb}
\usepackage{eucal}
\newtheorem{cor}[subsection]{Corollary}
\newtheorem{lem}[subsection]{Lemma}
\newtheorem{prop}[subsection]{Proposition}

\newtheorem{thm}[subsection]{Theorem}

\theoremstyle{remark}

\newcommand{\thmref}[1]{Theorem~\ref{#1}}
\newcommand{\secref}[1]{Sect.~\ref{#1}}
\newcommand{\lemref}[1]{Lemma~\ref{#1}}
\newcommand{\propref}[1]{Proposition~\ref{#1}}
\newcommand{\corref}[1]{Corollary~\ref{#1}}

\newcommand{\nc}{\newcommand}
\newcommand{\renc}{\renewcommand}
\nc\on{\operatorname}\nc\ssec{\subsection}
\nc\Ag{{\mathcal A}_G}\nc\ag{{\mathfrak a}_G}
\nc\Og{{\mathcal O}_{\check G}}
\nc\GG{\mathbb G}
\nc\ZZ{\mathbb Z}
\nc\QQ{\mathbb Q}
\nc\CC{\mathbb C}
\nc\C{\mathfrak C}
\nc\Af{\mathfrak A}
\nc\R{{\mathcal R}}
\nc\rr{\mathfrak r}
\nc\B{{\mathcal B}}
\nc\agd{\overset{\bullet}{\mathfrak a}_G}
\renc\a{\mathfrak a}
\nc\m{\mathfrak m}
\nc{\Cat}{{\mathfrak{C}(\CA,\CO)}}
\nc{\CatO}{{\mathfrak{C}(\CO,\CO)}}
\nc{\CatA}{{\mathfrak{C}(\CA,\CA)}}
\nc{\Catg}{{\mathfrak{C}(\CA_G,\CO_{\check{G}})}}
\nc{\Catsc}{{\mathfrak{C}(\CA_G,\CO_{\check G_{sc}})}}
\nc{\Catgd}{\overset{\bullet}{\mathfrak{C}}(\CA_G,\CO_{\check{G}})}
\nc{\bfHiggs}{{\mathbf H}{\mathbf i}{\mathbf g}{\mathbf g}{\mathbf s}}
\nc{\bfCam}{{\mathbf C}{\mathbf a}{\mathbf m}}
\nc{\bfBun}{{\mathbf B}{\mathbf u}{\mathbf n}}
\nc\F{{\mathcal F}}\nc{\CB}{{\mathcal B}}
\nc{\CO}{{\mathcal O}}
\nc{\Spec}{{\operatorname{Spec}}}
\nc{\CA}{{\mathcal A}}
\renewcommand{\L}{{\mathbf L}}
\nc{\W}{{\mathbf W}}
\nc{\CM}{{\mathcal M}}
\nc{\1}{{\mathbf 1}}
\nc{\CL}{{\mathcal L}}
\nc{\U}{{\mathbf U}}
\nc{\Bl}{{\mathsf {Bl}}}
\nc{\tu}{\overset{\bullet}{\mathfrak u}}
\nc{\TM}{\overset{\bullet}{M}}
\nc{\TL}{\overset{\bullet}{L}}
\renc{\lg}{\check \g}\nc{\Fr}{\on{Fr}}
\nc{\Gr}{\on{Gr}}
\nc{\Fl}{\on{Fl}}
\nc{\Ad}{\on{Ad}}
\nc\PG{\on{P}_{G[[t]]}(\Gr)}
\nc\PF{\on{P}_{G[[t]]}(\Fl)}
\nc\PFt{\widetilde{\on{P}}_{G[[t]]}(\Fl)}
\nc\PFti{\widetilde{\underline{\on{P}}}_{G[[t]]}(\Fl)}
\nc{\Res}{\operatorname{Res}}
\nc{\g}{{\mathfrak g}}
\nc\Iw{\on{Iw}}
\renewcommand{\u}{{\mathfrak u}}
\nc{\Ind}{\operatorname{Ind}}
\nc{\Ker}{\operatorname{Ker}}
\nc{\im}{\operatorname{Im}}
\nc{\Coker}{\operatorname{Coker}}
\nc\dirlim{\underset{\map}{\on{lim}}}
\nc\Flt{\widetilde{\Fl}}
\nc\DFlc{\on{D}^{G[[t]]}_c(\Flt)}
\nc\gh{\widehat\g}\nc\KMc{\gh^{G[[t]]}_c\mod_0}
\nc{\Tor}{\operatorname{Tor}_{\bullet}}
\nc{\mapr}{\operatorname{Tor}}
\nc{\Hom}{\operatorname{Hom}}

\renewcommand{\mod}{\operatorname{-mod}}
\nc{\mmod}{-\underline{\on{mod}}}
\nc{\comod}{\operatorname{-comod}}

\nc{\Barb}{\operatorname{Bar}^{\bullet}}
\nc{\ten}{{\otimes}}
\nc{\map}{\longrightarrow}

\begin{document}

\author{Sergey Arkhipov and Dennis Gaitsgory}

\title[Modules over the small quantum group]
{Another realization of the category of modules over the small quantum group}

\address{S.A.: Ind. Moscow Univ., 11 Bolshoj Vlasjevskij per., 
Moscow, 121002 Russia; \newline
D.G.: Dept. of Math., Harvard University, Cambridge MA 02138}

\email{hippie@mcmme.ru, gaitsgde@math.harvard.edu}

\date{}

\maketitle
\section*{Introduction}
\ssec{}
Let $\g$ be a semi-simple Lie algebra. Given a root of unity
(cf. \secref{biggroup}), one can consider two remarkable algebras,
$\U_\ell$ and $\u_\ell$, called the {\it big} and the {\it small} quantum
group, respectively. Let $\U_\ell\mod$ and $\u_\ell\mod$
denote the corresponding categories of modules.
It is explained in \cite{Lu4} and \cite{AJS} that the
former is an analog in characteristic $0$ of the category of algebraic
representations of the corresponding group $G$ over a field of positive
characteristic, and the latter is an analog of the category of
representations of its first Frobenius kernel.

It is  a fact of
crucial importance, that although $\U_\ell$ is introduced as an
algebra defined by an explicit set of generators and relations,
the category $\U_\ell\mod$ (or, rather, its regular block,
cf. \secref{reg}) can be described in purely geometric terms, as
perverse sheaves on the (enhanced) affine flag variety $\Flt$,
cf. \secref{KLKT}. This is obtained by combining the
Kazhdan-Lusztig equivalence between quantum groups and
affine algebras and the Kashiwara-Tanisaki localization of modules
over the affine algebra on $\Flt$. This paper is a first step in the
project of finding a geometric realization of the category $\u_\ell\mod$.
We should say right away that one such realization already exists, and is
a subject of \cite{BFS}. However, we would like to investigate 
other directions. 

We were motivated by a set of conjectures proposed by B.~Feigin, 
E.~Frenkel and G.~Lusztig, which, on the one hand,
tie the category $\u_\ell\mod$ to the (still hypothetical) category of
perverse sheaves on the semi-infinite flag variety
(cf. \cite{FM}, \cite{FFKM}), and on the other hand,
relate the latter to the category of modules over the affine algebra at
the critical level. 

Since we already know the geometric interpretation for modules 
over the big quantum
group, it is a natural idea to first express $\u_\ell\mod$ entirely
in terms of $\U_\ell\mod$. This is exactly what we do in this paper.

\ssec{} According to \cite{lubook}, there is a functor $\Fr^*$
from the category of finite-dimensional representations of the Langlands
dual group to $\U_\ell\mod$. In particular, we obtain a
bi-functor: $\check G\mod\times \U_\ell\mod\to\U_\ell\mod$: $V,M\to
\Fr^*(V)\ten M$. We introduce the category $\Catg$ to have as objects
$\U_\ell$-modules $M$, which satisfy the {\it Hecke eigen-condition}, in
the sense of \cite{BD}.

In other words, an object of $\Catg$ consists
of $M\in\U_\ell\mod$ and a collection of maps $\alpha_V: \Fr^*(V)\ten M
\to \underline{V}\ten M$, where $\underline{V}$ is the vector
space underlying the representation $V$. The main result of this paper
is \thmref{main}, which states that there is a natural equivalence between
$\Catg$ and $\u_\ell\mod$.

As the reader will notice, the proof of \thmref{main} is extremely simple. However,
it allows one to give the desired description of the regular block
$\u_\ell\mod_0$ of the category of $\u_\ell$-modules in terms of perverse
sheaves on the enhanced affine flag variety satisfying the Hecke
eigen-condition, cf. \secref{geominter}.

In a future publication,
we will explain how \thmref{geominter} can be used to define a
functor from $\u_\ell\mod_0$ to the category of perverse
sheaves on the semi-infinite flag variety and to other interesting
categories that arise in representation theory. In particular,
$\u_\ell\mod_0$ obtains an interpretation in terms of the geometric
Langlands correspondence: it can be thought of as a categorical
counterpart of the space of Iwahori-invariant vectors in a spherical
representation.

In another direction, \thmref{main} has as a
consequence the theorem that $\u_\ell\mod$ is equivalent to the category
of $G[[t]]$-integrable representations of the chiral Hecke algebra,
introduced by Beilinson and Drinfeld. (We do not state this theorem
explicitly, because the definition of the chiral Hecke algebra
is still unavailabale in the published literature.)

\ssec{}
Let us briefly describe the contents of the paper.

In Sect. 1 we recall the basic definitions concerning quantum groups.

In Sect. 2 we state our main theorem and its generalization for pairs
of bi-algebras $(\CA,\a)$.

In Sect. 3 we prove \thmref{main}
in the general setting.

In Sect. 4 we discuss several categorical
interpretations of \thmref{main} and, in particular, its variant
that concerns the graded version $\tu_\ell$ of $\u_\ell$.

In Sect. 5
we discuss the relation between the block decompositions of $\U_\ell$
and $\u_\ell$.

Finally, in Sect. 6 we prove \thmref{geominter}, which provides a
geometric interpretation for the category $\u_\ell\mod_0$.

\medskip

In this paper we consider  quantum groups at a root of unity of an even order,
in order to be able to apply the Kazhdan-Lusztig equivalence. 
However, the main result i.e. \thmref{main} holds and can be proved
in exactly the same way in the case of a root of unity of an odd order, with the
difference that in the definition of the quantum Frobenius, the Langlands
dual group $\check G$ must be replaced by $G$.

\ssec{Acknowledgments}
The main idea of this paper, i.e. \thmref{main}, occurred to us
after a series of conversations with B.~Feigin, M.~Finkelberg and A.~Braverman,
to whom we would like to express our gratitude.

In addition, we would like mention that \thmref{main},
was independently and almost simultaneously obtained by B.~Feigin and E.~Frenkel.

\section{Quantum groups} \ssec{Root data}

Let $G$ be a semi-simple simply-connected group. Let $T$ be the Cartan
group of $G$ and let $(I,X,Y)$ be the corresponding root data,
where $I$ is the set of vertices of the Dynkin diagram, $X$ is the set of
characters $T\map \GG_m$ (i.e. the weight lattice of $G$) and $Y$ is the set
of co-characters $\GG_m\map T$ (i.e. the coroot lattice of $G$).
We will denote by $\langle\, ,\rangle$ the canonical pairing
$Y\times X\map \ZZ$.
For every $i\in I$, $\alpha_i\in X$ (resp., $\check\alpha_i\in Y$)
will denote the corresponding simple root (resp., coroot);
for $i,j\in I$ we will denote by $a_{i,j}$ the corresponding
entry of the Cartan matrix, i.e. $a_{i,j}=\langle \check\alpha_i,
\alpha_j\rangle$.

Let $(\cdot,\cdot):X\times X\map\QQ$ be the {\it canonical} inner form.
In other words, $||\alpha_i||^2=2d_i$, where $d_i\in\{1,2,3\}$
is the minimal set of integers such that the matrix
$(\alpha_i,\alpha_j):=d_i\cdot a_{i,j}$ is symmetric.

\ssec{The big quantum group.} \label{biggroup}

Given the root data $(I,Y,X)$ Drinfeld and Jimbo constructed a Hopf
algebra $\U_v$ over the field $\CC(v)$ of rational functions in $v$.
Namely, $\U_v$ has as generators the elements
\footnote{We are using a slightly non-standard version of $\U$, in which the toric
part coincides with the group-algebra of the classical torus $T$}
$E_i, F_i,\ i\in I$, $K_t, t\in T$ and the relations are:

\begin{align*}
& K_{t_1}\cdot K_{t_2}=K_{t_1\cdot t_2}, \\
& K_t\cdot E_i\cdot K_t^{-1}=\alpha_i(t)\cdot E_i,\,\,\,K_t
\cdot F_i\cdot K_t^{-1}=\alpha_i(t^{-1})\cdot F_i  \\
& E_i\cdot F_j-F_j\cdot E_i=\delta_{i,j}\cdot\frac {K_i-K_i^{-1}}{v^{d_i}-
v^{-d_i}}, \text{ where } K_i=K_{d_i\cdot\check\alpha_i(v)}, \\
&\sum_{r+s=1-a_{ij}}(-1)^s{\left[{\begin{array}{c}1-a_{ij}\\ s\end{array}}
\right]}_{d_i}
E_i^r\cdot E_j\cdot E_i^s=0
\text{ if } i\ne j,\\
&\sum_{r+s=1-a_{ij}}(-1)^s{\left[{\begin{array}{c}1-a_{ij}\\ s\end{array}}
\right]}_{d_i}
F_i^r\cdot F_j\cdot F_i^s=0
\text{ if } i\ne j, \text{ where } \\
&\left[{\begin{array}{c}m\\ t\end{array}}\right]_d:=\prod_{s=1}^t{\frac
{v^{d\cdot (m-s+1)}-v^{-d\cdot (m-s+1)}} {v^{d\cdot s}-v^{-d\cdot s}}} \text{ for }
m\in \ZZ.
\end{align*}

The co-product is given by the formulae:
\begin{align*}
&\Delta(E_i)=E_i\otimes 1+K_i\ten E_i,\\
&\Delta(F_i)=F_i\otimes K_i^{-1}+1\ten F_i,\\
&\Delta(K_t)=K_t\ten K_t,
\end{align*}
and the co-unit and antipode maps are
\begin{align*}
&\epsilon(E_i)=\epsilon(F_i)=0,\;\;\;\;\epsilon(K_t)=1,\\
&\tau(K_t)=K_{t^{-1}},\;\;\;\; \tau(E_i)= -K_i^{-1}
\cdot E_i,\;\;\;\; \tau(F_i)=-F_i\cdot K_i.
\end{align*}

\smallskip

Let now $\ell$ be a sufficiently large even natural number, which divides all the
$d_i$'s. We set $\ell_i=\ell/d_i$ and let us fix a primitive $2\ell$-th
root of unity $\zeta$. Let $R\subset \CC(v)$ denote the localization
of the algebra $\CC[v,v^{-1}]$ at the ideal corresponding to $v-\zeta$.

In his book \cite{lubook}, G.~Lusztig defined an $R$-lattice $\U_R$ inside $\U_v$. Namely,
$\U_R$ is an $R$-subalgebra of $\U_v$ generated by $E_i,F_i,K_t$ and the following
additional elements:

\begin{align*}
& E_i^{(\ell_i)}:=\frac{E_i^{\ell_i}}{[\ell_i]_{d_i}!},\,\,
F_i^{(\ell_i)}:=\frac{F_i^{\ell_i}}{[\ell_i]_{d_i}!},
\text{  where   } [m]_d!=\prod_{s=1}^m {\frac{v^{d\cdot s}-v^{-d\cdot s}}{v^d-v^{-d}}}, \\
& \text{ and  } \left[{\begin{array}{c}K_i; m\\ t\end{array}}\right]_{d_i}:=\prod_{s=1}^t{\frac
{K_i\cdot v^{d_i\cdot (m-s+1)}-K_i^{-1}\cdot v^{-d_i\cdot (m-s+1)}} {v^{d_i\cdot s}-v^{-d_i\cdot s}}},
\text{ for } m\in \ZZ.
\end{align*}

It is shown in {\it loc.cit.} that $\U_R$ is a Hopf subalgebra of $\U_v$. Finally, following Lusztig we set
$\U_\ell$ to be the reduction of $\U_R$ modulo the ideal $(v-\zeta)\subset R$. By construction, $\U_\ell$
is a Hopf algebra over $\CC$.

\smallskip

The main object of study of this paper is not so much the
algebra $\U_\ell$ itself, but rather certain categories of
its representations. We introduce the category $\U_\ell\mod$
as follows: its objects are finite-dimensional representations
$M$ of $\U_\ell$, for which the action of the $K_t$'s comes from
an algebraic action of the torus $T$ on $M$, and such that 
for $\lambda\in X$, the action of 
$\left[{\begin{array}{c}K_i; m\\ t\end{array}}\right]_{d_i}$ on the
subspace of $M$ of weight $\lambda\in X$ is given by the scalar
$\left[{\begin{array}{c}\langle \check \alpha_i, \lambda\rangle +m \\ t\end{array}}\right]_{d_i}$.
(Note that the elements $\left[{\begin{array}{c}m\\ t\end{array}}\right]_{d_i}\in \CC(v)$
all belong to $R$, and hence they are well-defined in $\CC=R/(v-\zeta)$.)

The category $\U_\ell\mod$ is a monoidal
category endowed with a forgetful functor to the category of finite-dimensional $\CC$-vector spaces. Hence,
there exists a Hopf algebra, such that the category $\U_\ell\mod$
is equivalent to the category of finite-dimensional co-modules over it. We will denote
this Hopf algebra by $\Ag$.

One should think of $\Ag$ as of a quantization
of the algebra of regular functions on the group $G$. It is known that
$\Ag$ is finitely generated as an associative algebra. Moreover, we will
see that $\Ag$ contains a large commutative subalgebra,
over which it is finitely generated as a module.

\ssec{Quantum Frobenius homomorphism} \label{quantum Frob}

Let $(I,X^*,Y^*)$ be the Langlands dual root data. In other words,
$X^*:=Y$ and $Y^*:=X$ are the weight and the coweight lattices of
the Langlands dual torus $\check T$. The corresponding semi-simple group
$\check G$ is by definition of the adjoint type. Let $\lg$
denote the Lie algebra of $\check G$. Let $\check G\mod$ denote the
category of finite-dimensional $\check G$-modules and let $\CO_{\check G}$
be the algebra of functions on $\check G$. We will denote by $\U(\lg)$ the usual
universal enveloping algebra of $\lg$.

\medskip

The canonical inner form $(\cdot,\cdot)$ on $X$ gives rise to the inner
form on $Y$, which is not necessarily integral-valued, since
$||\check\alpha_i||=\frac{2}{d_i}$. However, if we multiply the
latter by $\ell$, we obtain an integral valued form
$(\cdot,\cdot)_\ell: Y\times Y\map\ZZ$.
By construction,
$$||\check\alpha_i||_\ell^2=2\cdot \ell_i \text{ and }
(\check\mu,\check\alpha_i)_\ell=\ell_i\cdot \langle \check\mu,
\alpha_i\rangle.$$
Using the pairing $(\cdot,\cdot)_\ell$
we obtain the map $\phi:Y\map X$ given by $\check\mu\mapsto
(\check\mu,\cdot)_\ell$ and the map $\phi_T:T\map \check T$.

\smallskip

Following Lusztig (\cite{lubook}, Theorem 35.1.9) one defines the
{\it quantum Frobenius morphism}. For us, this will be a functor
$$\Fr^*:\check G\mod\map\U_\ell\mod,$$
constructed as follows:

Starting with a $\check G$-module $V$, we define a $\U_\ell$-action on it
by letting the torus $T$ act via
$$T\overset{\phi_T}\to \check T\hookrightarrow \check G,$$
which defines the action of the $K_t$'s and the 
$\left[{\begin{array}{c}K_i; m\\ t\end{array}}\right]_{d_i}$'s.

The generators $E_i, F_i$ will act by $0$, and $E_i^{(\ell_i)}, F_i^{(\ell_i)}$ will act as the corresponding
Chevalley generators $e_i$ and $f_i$ of $\U(\lg)$.

It is essentially a theorem of Lusztig, (\cite{lubook}, Theorem 35.1.9) that the above formulae
indeed define an action of $\U_\ell$ on $V$. Moreover, from {\it loc.cit.} it follows that
the functor $\Fr^*$ preserves the tensor structure and is full. Hence, we obtain an injective
homomorphism of Hopf algebras $\phi_G:\Og\map \Ag$.

\medskip

Let $\check G_{sc}$ be the simply-connected cover of the
group $\check G$ and let $X^*_{sc},Y^*_{sc}$ and $\check T_{sc}$
be the corresponding objects for $\check G_{sc}$. In particular,
$Y$ identifies with the coroot lattice inside the coweight lattice 
$X^*_{sc}$, and $Y^*_{sc}=\on{Span}(\alpha_i)$.

Since $\phi(\check\alpha_i)=\ell_i\cdot \alpha_i$, we obtain
that the map $\phi$ gives rise to a map $\phi_{sc}:X^*_{sc}\to X$.

Therefore, we have a map $\phi_{T,sc}:T\to \check T_{sc}$
and the functor $\on{Fr}^*:\check G\mod\to \U_\ell\mod$ can
be extended to a tensor functor 
$\on{Fr}^*_{sc}:\check G_{sc}\mod\to \U_\ell\mod$ by the same
formula.

\ssec{The small quantum group}   \label{small quantum group}

Following Lusztig, we first define the graded version of the small
quantum group, denoted $\tu_\ell$.

By definition, this is a sub-algebra of $\U_\ell$ generated by 
$E_i, F_i, i\in I$ and all the $K_t$'s.

From the formula for coproduct of the above generators, it follows that $\tu_\ell$
is in fact a Hopf subalgebra of $\U_\ell$.

We define the category $\tu_\ell\mod$ to consist of all finite-dimensional
$\u_\ell$-modules $M$, on which the action of the $K_t$'s comes from an 
algebraic action of $T$ on $M$.

The restriction functor $\U_\ell\mod\to\tu_\ell\mod$ corresponds
to a map of Hopf algebras $\Ag\to\agd$.



\medskip

Finally, we are ready to introduce our main object of study--
the {\it small quantum group}, $\u_\ell$. One would want it to be
a Hopf subalgebra of $\U_\ell$, universal with the 
property that it acts trivially on representations of the form $\Fr^*(V)$.

When one works with a root of unity of an odd order, the corresponding
subalgebra is just generated by $K_i$, $E_i$ and $F_i$, $i\in I$. However,
in the case of a root of unity of an even order considered in the present
paper, it appears that a {\it Hopf subalgebra} with such properties does
not exist.

In our definition, $\u_\ell$ will be just an associative subalgebra
of $\U_\ell$, generated by $K_iE_i,\, F_i, i\in I$ and $K_t$ for
$K_t\in\on{ker}(\phi_T)$. It is easy to see that $\u_\ell$
is finite-dimensional. 

We define the category $\u_\ell\mod$ to have as objects all 
finite-dimensional $\u_\ell$-modules.

By construction, we have a restriction functor 
$\Res:\U_\ell\mod\to \u_\ell\mod$. It corresponds to a 
homomorphism of co-algebras $\Ag\to \ag$. 

\medskip

Note that although the co-product on $\U_\ell$ does not preserve
$\u_\ell$, it maps it to $\u_\ell\otimes \U_\ell$. This means that
$\ag$ has a structure of a right $\Ag$-module. In categorical terms,
we have a well-defined functor $(M\in\u_\ell\mod,\, 
N\in \U_\ell\mod)\mapsto
M\otimes \Res(N)\in \u_\ell\mod$.

\medskip

\noindent{\it Remark.}
As we shall see later, although $\ag$ is not a Hopf algebra,
the category $\u_\ell\mod$ will be in fact a monoidal category.
The ``paradox'' is explained as follows: the tautological forgetful
functor $\u_\ell\mod\to \{\text{Vector spaces}\}$ cannot be made into
a tensor functor so that the composition
$$\U_\ell\mod\overset{\Res}\to \u_\ell\mod\to \{\text{Vector spaces}\}$$
is the standard fiber functor on $\U_\ell\mod$.

\medskip

Consider now the restriction of the quantum Frobenius to $\u_\ell$, i.e.
the composition
$$\check G\mod\overset{\Fr^*}\map \U_\ell\mod\overset{\on{Res}}
\map \u_\ell\mod.$$

From the formula for $\Fr^*$ it is easy to see that it factors through
the forgetful functor $V\mapsto\underline{V}$ from 
$\check G\mod$ to vector spaces, i.e.
$$\check G\mod\to \{\text{Vector spaces}\} \overset{\on{co-unit}}\to \u_\ell\mod.$$
Moreover, we have the following assertion (\cite{lubook}, Theorem 35.1.9):

\begin{prop} \label{factorization}
Let $M$ be an object of $\U_\ell\mod$. Then

\smallskip

\noindent{\em (1)} The subspace of $\u_\ell$-invariants
$M^{\u_\ell}\subset M$ (i.e. $m\in M^{\u_\ell}$ if 
$u\cdot m=\epsilon(u)\cdot m$ for $u\in \u_\ell$)
is $\U_\ell$-stable.

\smallskip

\noindent{\em (2)} If the $\u_\ell$-action on $M$ is trivial, 
there exists a (unique up to a unique
isomorphism) $\check G$-module $V$ such that $M\simeq \Fr^*(V)$.
\end{prop}

For completeness, let us sketch the proof of the second part 
of this proposition.

\begin{proof}

From the short exact sequence $1\to \on{ker}(\phi_T)\to T\to \check T\to 1$
we obtain the $T$-action on $M$ comes from a $\check T$-action.

In particular, the element $\left[{\begin{array}{c}K_i; 0\\ \ell_i\end{array}}\right]_{d_i}$
acts on $M$ as the Lie algebra element $h_i\in \check \g$.

We define the action of $e_i$ and $f_i$ as $E_i^{(\ell_i)}$ and $F_i^{(\ell_i)}$, 
respectively, and we need just to check that the relation $[e_i,f_i]=h_i$ holds.
But this follows from the formula
$$[E_i^{(\ell_i)},F_i^{(\ell_i)}]=\underset{0\leq k<\ell_i}\Sigma\,\,
(\frac{1}{[k]_{d_i}!})^2\cdot (E_i)^k\cdot \left[{\begin{array}{c}K_i; 2k\\ 
\ell_i-k\end{array}}\right]_{d_i}
\cdot (F_i)^k,$$
and all the terms but $\left[{\begin{array}{c}K_i; 0\\ \ell_i\end{array}}\right]_{d_i}$
belong to the two-sided ideal generated by the $E_i$'s and the $F_i$'s.

\end{proof}

\medskip

In a similar fashion one defines the ``simply-connected'' version of
$\u_\ell$, which we will denote by $\u_{\ell,sc}$. By definition,
this is an associative subalgebra of $\u_\ell$
generated by $K_iE_i,\, F_i, i\in I$ and $K_t$ for
$K_t\in\on{ker}(\phi_{T,sc})$.

The category $\u_{\ell,sc}\mod$ and the co-algebra $\ag{}_{sc}$ are defined
in a similar way. The analog of \propref{factorization}  
above holds for $\u_\ell$ replaced by $\u_{\ell,sc}$ and $\check G$ replaced by 
$\check G_{sc}$, respectively.

\section{The main result}

\ssec{The category $\Catg$}

We now come to the definition crucial for this paper. To avoid
redundant repetitions, we will work with $\u_\ell$ (resp., $\check G$),
while the case of $\u_{\ell,sc}$ (resp., $\check G_{sc}$) can be treated
similarly.
 
Let us consider the ind-completions of the categories $\U_\ell\mod$,
$\u_\ell\mod$ and $\check G\mod$. Each of these categories consists of
all co-modules over the corresponding co-algebra, i.e. $\Ag$,
$\ag$ or $\Og$, respectively.

We define the category $\Catg$ to have as objects
vector spaces $M$ endowed with an action of the {\em algebra} $\Og$ and
with a co-action of the {\em co-algebra} $\Ag$ compatible in the
following natural way:
$$\on{co-ac}(f\cdot m)=\Delta(f)\cdot(\on{co-ac}(m)).$$
Here $f\in\Og$, $m\in M$, $\on{co-ac}:\ M\map\Ag\ten M$ denotes the
co-action map, the element $\Delta(f)$ belongs to
$\Og\ten\Og\subset\Ag\ten\Ag$ and acts on $\Ag\ten M$. Morphisms in
this category are the ones preserving both the action and the co-action.

\medskip

In other words, we need that the action map $\Og\otimes M\to M$
is a map of $\Ag$-comodules, or equivalently, that the co-action map
$M\to \Ag\ten M$ is the map of $\Og$-modules.

An example of an object of $\Catg$ is $M=\Og$, with the natural
$\Ag$-coaction (coming from the fact that $\Og$ is a Hopf subalgebra
in $\Ag$) and the $\Og$-action. Another basic example is $M=\Ag$.

\ssec{A reformulation.} \label{reform}
Here is a more ``geometric'' way to formulate this definition.
We claim that the category $\Catg$ is equivalent to the category of
pairs 
$$(M\in \Ag\comod,\{\alpha_V,\,\,\forall\, V\in \check G\mod\}),$$
where each $\alpha_V$ is a map of $\Ag$-comodules
(i.e., of $\U_\ell$-modules)
$$\alpha_V: \Fr^*(V)\otimes M\to \underline{V}\otimes M$$
(recall that for $V\in \check G\mod$, the notation 
$\underline{V}$ stands for the underlying vector space),
such that
\begin{itemize}

\item For $V=\CC$, $\alpha_V: M\to M$ is the identity map.

\item For a map $V_1\to V_2$, the diagram
$$\CD\Fr^*(V_1)\otimes M  @>{\alpha_{V_1}}>>  \underline{V_1}\otimes M  \\
@VVV        @VVV    \\
\Fr^*(V_2)\otimes M  @>{\alpha_{V_2}}>>  \underline{V_2}\otimes M\endCD$$ commutes.

\item
A compatibility with tensor products holds in the sense that the map
$$\Fr^*(V_1)\otimes \Fr^*(V_2) \otimes M\to \Fr^*(V_1\ten V_2)
\otimes M\overset{\alpha_{V_1\otimes V_2}}
\longrightarrow\underline{V_1\otimes V_2}
\otimes M\to \underline{V_2}\otimes \underline{V_1}\otimes M$$
equals
$$\Fr^*(V_1)\otimes \Fr^*(V_2)\otimes M\overset{\on{id}\otimes\alpha_{V_2}}
\longrightarrow \Fr^*(V_1)\otimes \underline{V_2}\otimes M\simeq
\underline{V_2}\otimes \Fr^*(V_1)\otimes M\overset{\on{id}\otimes\alpha_{V_1}}\longrightarrow
\underline{V_2}\otimes \underline{V_1}\otimes M.$$

\end{itemize}

Morphisms in this category between $(M,\alpha_V)$ and $(M',\alpha'_V)$
are $\U_\ell$-module maps $M\to M'$, such that each square
$$\CD\Fr^*(V)\otimes M   @>{\alpha_V}>> \underline{V}\otimes
M \\@VVV       @VVV        \\ \Fr^*(V)\otimes M'   @>{\alpha'_V}>>
\underline{V}\otimes M'\endCD$$
commutes.

\bigskip

Indeed, given $M$ as above we define the action of $\Og$ on it as the
composition map $$\Fr^*(\Og)\otimes M\overset{\alpha_{\Og}}
\longrightarrow \underline{\Og}
\otimes M\overset{\epsilon\ten\on{id}}\longrightarrow M,$$
where $\epsilon$ is the co-unit $f\mapsto f(1)$ in $\Og$.
Conversely, given $M\in \Catg$, the map $\alpha_V$ comes by
adjunction from the map
$$\Fr^*(V)\otimes \underline{V}^*\otimes M
\overset{\text{matr.coef.}
\otimes\on{id}}\longrightarrow\Fr^*(\Og)\otimes M\to M.$$
Let us make the following observation:

\begin{prop}
\label{isom}
For $(M,\alpha_V)\in\Catg$, the maps $\alpha_V$ are automatically
isomorphisms.

\end{prop}

\begin{proof}
Let $N$ be the kernel of the map $\Fr^*(V)\ten M\to \underline{V}\ten M$
and let $V^*\in \check G\mod$ be the dual of $V$. From
the axioms on the $\alpha_V$'s,
we obtain that the composition $$\Fr^*(V^*)\ten N\to \Fr^*(V^*\ten V)\ten
M\to \underline{V^*\ten V}\ten M\to M$$
is on the one hand zero, and on the other hand equals the natural
 map $\Fr^*(V^*)\ten N\to M$, which is a contradiction. The surjectivity
of $\alpha_V$ is proved in the same way.
\end{proof}

Our main result is the following theorem:

\begin{thm} \label{main}
The category $\Catg$ is naturally equivalent to the category of
$\ag$-comodules. Objects in $\Catg$, which are finitely generated over
$\Og$, correspond under this equivalence to finite-dimensional
$\ag$-comodules.
\end{thm}

This theorem has the following interesting corollary:

\begin{cor} \label{dual group action}
The Langlands dual group $\check G$ acts on the category
$\u_\ell\mod$ by endo-functors. In other words, 

\smallskip

\noindent {\em (1)}
For every $\gamma\in \check G$ there is a functor 
$T_\gamma:\u_\ell\mod\to \u_\ell\mod$.

\smallskip

\noindent {\em (2)}
For each pair $\gamma_1,\gamma_2\in\check G$
there is an isomorphism of functors
$T_{\gamma_1}\circ T_{\gamma_2}\Rightarrow T_{\gamma_1\cdot \gamma_2}$.

\smallskip

\noindent {\em (3)}
For each triple $\gamma_1,\gamma_2,\gamma_3$ the two
natural transformations 
$T_{\gamma_1}\circ T_{\gamma_2}\circ T_{\gamma_3}\Rightarrow 
T_{\gamma_1\cdot \gamma_2\cdot \gamma_3}$ coincide.

\end{cor}

\begin{proof}

Let us view $\u_\ell$-modules as objects of $\Catg$ via
\thmref{main}.

Given an object $(M,\alpha_V)\in \Catg$ and an element 
$\gamma\in \check G$ we define a new object $T_\gamma(M,\alpha_V)$
as follows:

The underlying $\U_\ell$-module is the same, i.e. $M$. However,
the corresponding morphism
$$\on{Fr}^*(V)\otimes M\to \underline{V}\otimes M$$
is the old $\alpha_V$ composed with 
$\underline{V}\otimes M\overset{\gamma\otimes\on{id}}\longrightarrow 
\underline{V}\otimes M$, where $\gamma\in \check G$ is viewed
as an automorphism of the vector space $\underline{V}$.

It is clear that in this way we indeed obtain an action
of $\check G$ on $\Catg$, and hence on $\ag$-comod, 
by endo-functors. 

\end{proof}

Another corollary of \thmref{main} is as follows:

\begin{cor}  
The category $\u_\ell\mod$ has a natural monoidal structure.
\end{cor}

\begin{proof}

Given two objects $(M,\alpha_V)$ and $(M',\alpha'_V)$ in $\Catg$
we have to define their tensor product $(M'',\alpha''_V)$
as a new object of $\Catg$.

Consider first their naive tensor product $M\otimes M'$ as a $\U_\ell$-module.
We claim that the algebra $\Og$ acts on it by endomorphisms. Indeed, to define
such an action, it is enough to define $\U_\ell$-module maps
$$\underline{V}\otimes (M\otimes M')\to \underline{V}\otimes (M\otimes M')$$
for every $V\in \check G\mod$, compatible with the tensor structure on 
$\check G\mod$ in the same sense as in the definition of $\Catg$.

The sought-for maps are defined as follows:
$$\underline{V}\otimes M\otimes M'\overset{\alpha_V}\simeq 
(\Fr^*(V)\otimes M)\otimes M'\simeq M\otimes (\Fr^*(V)\otimes M')
\overset{\alpha'_V}\simeq \underline{V}\otimes M\otimes M',$$
where the second arrow comes from the braiding on the category
$\U_\ell\mod$.

The $\U_\ell$-module $M''$ is defined as the fiber at $1\in \check G$
of $M\otimes M'$ viewed as a quasi-coherent sheaf on $\check G$.
It comes equipped with a data of $\alpha''$ by construction.

It is easy to see that the functor $(M,\alpha_V),(M',\alpha')\mapsto 
(M'',\alpha''_V)$ admits a natural associativity constraint, which makes
$\Catg$ into a monoidal category. Moreover, if both $M$ and $M'$ are
finitely generated as $\Og$-modules, then so is $M''$. Hence, this monoidal
structure preserves the sub-category of finite-dimensional $\ag$-comodules,
which is the same as $\u_\ell\mod$.

\end{proof}

\ssec{The general setting} \label{coalset}

It will be convenient to generalize our setting as follows.
Let $\CO$, $\CA$ be two Hopf algebras and let $\CO\to\CA$ be an embedding.

In addition, let $\a$ be a co-algebra and a right $\CA$-module,
and let $\CA\to\a$ be a surjection respecting both structures.

We impose the following conditions on our data:

\medskip

\begin{itemize}

\item[(i)]
The composition $\CO\to\CA\to\a$ factors as
$\CO\overset{\on{co-unit}}\longrightarrow
\CC\overset{\on{unit}}\longrightarrow \a$.

\medskip

\item[(ii)]
The inclusion $\CO\subset \CA^\a$ is an equality.
\footnote{For a 
co-algebra $\B$ co-acting on $M$, the notation $M^\B$ will mean ``invariants'',
i.e. $$M^\B=\Hom_{\B\comod}(\CC,M)=\on{Ker}(M\overset{\on{co-ac}-1
\otimes \on{id}}\map \B\ten M).$$
For example, for a $\u_\ell$-module
$M$ (which is the same as an $\ag$-comodule) $M^{\u_\ell}=M^{\ag}$.}

\medskip

\item[(iii)]
The inclusion $\m\cdot \CA\subset \on{Ker}(\CA\to\a)$
is an equality, where $\m$ is the augmentation ideal in
$\CO$.

\end{itemize}

In addition, we impose the following technical condition, that
one of the following two properties is satisfied (compare 
with \cite{Mo}, Sect. 3.4):

\medskip

\begin{itemize}

\item[(iv a)]
Either $\CA$ is faithfully-flat as an $\CO$-module,

\medskip

\item[(iv b)]
or the induction functor $\Ind:\a\comod\to\CA\comod$
(cf. \secref{induction functor} for the definition of $\Ind$)
is exact and faithful.

\end{itemize}

Of course, we will prove that our triple $(\Og,\Ag,\ag)$ satisfies
conditions (i-iv).

\medskip

We define  the category
$\Cat$ to have as objects vector spaces $M$ endowed with a left action of
the {\em algebra} $\CO$ and a left co-action of the {\em co-algebra}
$\CA$ which are compatible in the same sense as
in the definition of $\Catg$.
The following is a generalization of \thmref{main}:

\begin{thm} \label{gener}
The categories $\Cat$ and $\a\comod$ are naturally
equivalent.
\end{thm}

\ssec{The ``classical'' case}
The general \thmref{gener} models the following familiar situation.
Let $$1\to H'\to H''\to H\to 1$$
be a short exact sequence of linear algebraic groups.
Take
$\CO=\CO_H$, $\CA=\CO_{H''}$ and $\a=\CO_{H'}$. Conditions
(i)-(iii) obviously hold and we
claim, that the assertion of \thmref{gener} in this case is following
well-known phenomenon:

\medskip

First, by definition, the category $\Cat$ is naturally equivalent to
the category $\on{QCoh}^{H''}(H)$ of $H''$-equivariant
quasi-coherent sheaves on $H$. By taking the fiber of a sheaf
at $1\in H$ we  obtain a functor $$\on{QCoh}^{H''}
(H)\to \on{QCoh}^{H'}(\on{pt}),$$
which is known to be an equivalence of categories.
However, $$\on{QCoh}^{H'}(\on{pt})\simeq H'\mod\simeq\a\comod.$$
The proof of \thmref{gener} in the general case will be essentially
a translation of the above two-line proof into the language of
Hopf algebras.

\medskip

\noindent {\it Remark.} 
Suppose that in the setting of \thmref{gener},
$\CO$ is in fact commutative, i.e. $\CO$ is
the algebra of functions on an affine group-scheme
$\Gamma$. 

Then we have an analog of {dual group action}, that
$\Gamma$ acts on the category $\a\comod$ by endo-functors.
In the above example of $(\CA=\CO_{H''},\a=\CO_{H'})$, 
this action corresponds
to the natural map of $\Gamma=H$ to the group of {\it outer}
automorphisms of $H'$.

\section{Proof of the main theorem}

We will first prove the general \thmref{gener}. Then we
show that conditions (i)-(iv) are satisfied for 
$\Ag,\Og$ and $\ag$.

\ssec{The functor of (finite) induction}   \label{induction functor}

Now we proceed to the proof of \thmref{gener} in general.

Let us recall the definition of the (finite) induction 
functor $\a\comod\mapsto\CA\comod$.

Recall that if $M_1^r$ is a right co-module and $M_2^l$ is a left co-module
over a co-algebra $\a$, it makes sense to consider the vector space
$(M_1^r\otimes M_2^l)^\a$, equal by definition to the equalizer of the two maps
$$\Delta_1\otimes \on{id},\,\on{id}\otimes \Delta_2:
M_1^r\otimes M_2^l\to M_1^r\otimes\a\otimes M_2^l.$$

For $M\in\a\comod$ consider $\CA$ as a left $\CA$-co-module and a right $\a$-comodule,
and set $\Ind(M):=(\CA\otimes M)^\a$, which carries a left $\CA$-coaction by functoriality.

By construction, this functor is left exact, and it is the right adjoint of the natural 
restriction functor $\Res:\CA\comod\to \a\comod$.

\medskip

Note now that since $\CO=\CA^\a$, the action of $\CO$ on $\CA\otimes M$ 
by left multiplication maps the 
subspace $(\CA\otimes M)^\a$ to itself. Moreover, this action $\CO$-action on 
$\Ind(M)$ is compatible with the $\CA$-coaction.
Therefore, the functor $\Ind$ can be extended to a functor 
from $\a\comod$ to $\Cat$, which we will denote by $\underline\Ind$.

\medskip

Let us consider two examples. First, it is easy to see that
$\underline\Ind(\a)\simeq\CA$. Secondly, 
$\underline\Ind(\CC)\simeq \CO$, and more generally, for
$M$ is of the form $\Res(N)$ for $N\in \CA\comod$, we have:
$\underline\Ind(\Res(N))\simeq \CO\otimes N$, with the
diagonal $\CA$-coaction and the $\CO$-action on the first factor.

\ssec{The adjoint functor}
Now we will define a functor $\Cat\to \a\comod$.

Given an object $N\in\Cat$, consider the vector space
$\Psi(N):=\CC\underset{\CO}\otimes N$, where $\CO\to\CC$ is the co-unit map.

Since the $\CO$-action on $N$ commutes with the $\a$-coaction,
$\Psi(N)$ carries a natural co-action of $\a$.

Thus, we obtain a functor, denoted $\underline{\Res}:\Cat\map \a\comod$.
By construction, this functor is right exact. 

\medskip

By definition, for $\CO$ viewed as an object of $\Cat$,
$\underline{\Res}(\CO)\simeq \CC$. Property (iii) of \secref{coalset}
implies that $\underline{\Res}(\CA)\simeq \a$.
More generally, for objects of $\Cat$ of the form $\CO\otimes\Res(N)$,
we have: $\underline{\Res}(\CO\otimes \Res(N))\simeq \Res(N)$.

\begin{prop} \label{adj}
The functor $\underline{\Res}$ is the left adjoint to
$\underline{\Ind}$.

\end{prop}

\begin{proof}
We need to construct adjunction maps
$$\underline{\Res}\circ \underline{\Ind}(M)
\to M \text{ and }N\to \underline{\Ind}\circ \underline{\Res}(N)$$
for $M$ and $N$ in $\a\mod$ and $\Cat$, respectively.

Let $M$ be as above. Consider the composition
$(\CA\otimes M)^\a\hookrightarrow
\CA\otimes M\overset{\epsilon\otimes\on{id}}\longrightarrow M$.
By construction, this is a map of $\a$-comodules and it
obviously factors through
$$(\CA\otimes M)^\a\to \Psi((\CA\otimes M)^\a)\to M.$$
Therefore, we obtain a map $\underline{\Res}\circ \underline{\Ind}(M)\simeq
\Psi((\CA\otimes M)^\a)\to M$.

\medskip

For $N\in\Cat$, consider the map 
$$N\overset{\Delta}\to (\CA\otimes N)^\CA\hookrightarrow (\CA\otimes N)^\a\to
\Psi((\CA\otimes N)^\a).$$
This map respects the $\CA$-coaction and the $\CO$-action
by construction.

Thus, we obtain a map
$$N\to (\CA\otimes \Psi(N))^\a\simeq \underline{\Ind}\circ \underline{\Res}(N).$$

\end{proof}

Now we are ready to prove \thmref{gener}. We will give two proofs
corresponding to the two variants of condition (iv).

\ssec{Proof 1}

Let us first prove \thmref{gener} under the assumption that $(\CO,\CA,\a)$
satisfies condition (iv b) of \secref{coalset}, i.e.
that the induction functor $\Ind$ is exact and faithful.

We claim that the adjunction map
$$N\to \underline\Ind(\underline\Res(N))$$
is an isomorphism for any $N\in \Cat$.

Since the functor $\underline\Res$ is right-exact and $\Ind$ is exact,
the composition $\underline\Ind\circ\underline\Res$ is also right exact.
Hence, it suffices to show that for any $N$ as above there exists another
object $N'\in \Cat$ with a surjection $N'\twoheadrightarrow N$, for which
the map $N'\to \underline\Ind(\underline\Res(N'))$ is an isomorphism.

We set $N'=\CO\otimes N$, where $\CO$ {\it acts} on $\CO\otimes N$
via
$$a'\cdot (a\otimes n)\mapsto a'\cdot a\otimes n,$$
and the $\CA$-coaction is the diagonal one.
The map $\CO\otimes N\to N$ is given by the original $\CO$-action
on $N$.

Now, $\underline\Res(N')\simeq \Res(N)$, and 
$\underline\Ind(\underline\Res(N'))\simeq \CO\otimes N$, such that
the above adjunction map for $N'$ becomes the identity map 
on $\CO\otimes N$.

\medskip

Thus, to prove \thmref{gener}, it suffices to check that
the other adjunction map $\underline\Res(\underline\Ind(M))\to M$
is an isomorphism for any $M\in\a\comod$.
However, since the functor $\Ind$ (and hence $\underline\Ind$)
is faithful, it suffice to check that
$$\underline\Ind(\underline\Res(\underline\Ind(M)))\to 
\underline\Ind(M)$$
is an isomorphism. However, we know that the composition
$$\underline\Ind(M)\to
\underline\Ind(\underline\Res(\underline\Ind(M)))\to 
\underline\Ind(M)$$
is the identity map on $\underline\Ind(M)$, and the first arrow
is an isomorphism by what we have proved above. Hence, the second
arrow is an isomorphism as well.

\ssec{Proof 2}

Now let us prove \thmref{gener} under the assumption that $(\CO,\CA,\a)$
satisfies condition (iv a) of \secref{coalset}.

\begin{prop} \label{faithful flatness}
The functor $\underline{\Res}:\Cat\to\a\comod$
is exact and faithful.
\end{prop}

\begin{proof}

For an object $N\in \Cat$, consider the tensor product
$\CA\underset{\CO}\otimes N$. This is a left $\CA$-module
and a left $\CA$-comodule via the diagonal co-action.

Thus, we obtain a functor $\on{Coind}_\CO^\CA:\Cat\to\CatA$, which is exact and
faithful, since $\CA$ was assumed faithfully flat over $\CO$.

Now, the functor $N\mapsto \Psi(N)$ considered as a functor
from $\Cat$ to the category of vector spaces can be factored
as 
$$\Psi\simeq \Psi_{\CA}\circ \on{Coind}_\CO^\CA,$$
where $\Psi_{\CA}$ is the corresponding functor for $\CatA$.
Therefore, it suffices to show that $\Psi_{\CA}$ is exact and faithful.

However, the triple $(\CA,\CA,\CC)$ satisfies assumption (iii b), and we 
already know that $\Psi_{\CA}$ induces an equivalence between 
$\CatA$ and the category of vector spaces. In particular, $\Psi$
is exact and faithful.

\end{proof}

The rest of the proof proceeds very much in the same way as 
Proof 1 above.

\medskip

First, we claim that \propref{faithful flatness} above implies that the adjunction morphism
$$\underline{\Res}\circ\underline{\Ind}(M)\to M$$
is an isomorphism for every $M\in \a\comod$.
Indeed, every object in $\a\comod$ can be embedded
into a direct sum of several copies of $\a$, viewed as a co-module
over itself.
Hence, every $M\in \a\comod$ admits a resolution of the form:

$$M\to \a\otimes W_0\to \a\otimes W_1\to ...,$$
where $W_i$ are some vector spaces. Since the composition
$\underline{\Res}\circ\underline{\Ind}(M)$ is left-exact, it
is enough to prove that
$\underline{\Res}\circ\underline{\Ind}(\a)\to \a$ is an isomorphism.
However, this is obvious, since this map is the composition
$\underline{\Res}\circ\underline{\Ind}(\a)\simeq 
\underline{\Res}(\CA)\simeq \a$.

Thus, it remains to show that the adjunction map
$N\to \underline{\Ind}\circ \underline{\Res}(N)$
is an isomorphism. However, since the functor $\underline{\Res}$ is faithful,
it is enough to show that
$$\underline{\Res}(N)\to \underline{\Res}\circ \underline{\Ind}
\circ \underline{\Res}(N)$$
is an isomorphism. But we already know that
$\underline{\Res}\circ \underline{\Ind}\circ \underline{\Res}(N)\to
\underline{\Res}(N)$
is an isomorphism and the composition
$$\underline{\Res}(N)\to \underline{\Res}\circ
\underline{\Ind}\circ \underline{\Res}(N)\to\underline{\Res}(N)$$
is the identity map.

\medskip

\noindent{\it Remark}
Note that \thmref{gener} implies that under the assumption that
$\CA$ is faithfully-flat over $\CO$, the induction functor $\Ind:\a\comod\to\CA\comod$
is automatically exact. I.e., condition (iv a) in fact implies condition (iv b).

\ssec{Finiteness properties of $\CA$}
Thus, \thmref{gener} is proved. Let us now describe
the image of the category of finite-dimensional $\a$-comodules under our
equivalence of categories.

\begin{prop}
An $\a$-comodule $M$ is finite-dimensional if and only of
$\underline{\Ind}(M)$ is finitely generated as an $\CO$-module.

\end{prop}

\begin{proof}
One direction is clear: if $N\in\Cat$ is finite as an $\CO$-module,
then $\underline{\Res}(N)=\Psi(N)$ is finite-dimensional. 

Conversely, assume that $M$ is finite dimensional, and let
$M'\subset \underline\Ind(M)$ be a finite-dimensional
$\CA$-subcomodule, which surjects onto $M$ under
$$\Res(\underline\Ind(M))\to \underline\Res(\underline\Ind(M))\simeq M.$$

Then the $\CO$-submodule $N'$ in $\underline\Ind(M)$ generated by $M'$
is stable under both the $\CO$-action and $\CA$-coaction, and
$\underline\Res(N')$ surjects onto $M$. Hence, $N'=\underline\Ind(M)$.

\end{proof}

This proposition implies among the rest, that if $\a$ is
finite-dimensional, then
$\CA\simeq \underline{\Ind}(\a)$ is finitely generated as
a module over $\CO$. In particular, we obtain that in the quantum
group setting, $\Ag$ is a finite $\Og$-module.

\ssec{Verification of properties (i)-(iv) for quantum groups}

First, the map $\Og\to \Ag$ is injective because the quantum Frobenius 
homomorphism is surjective.

\medskip

To show that $\Ag\to\ag$ is surjective is equivalent to showing that every
object $M\in\u_\ell\mod$ appears as a sub-quotient of one of the form $\Res(N)$ for some
$N\in\U_\ell\mod$. We will prove a stronger assertion, namely, that any $M$ as above
is in fact a quotient of some $\Res(N)$:

\begin{prop}
For any $M\in\u_\ell\mod$, the canonical map
$\Res(\Ind(M))\to M$ is surjective.
\end{prop}

\begin{proof}

First, it is known (cf. \cite{APW1} Theorem 4.8 or \cite{APar} Proposition 3.15)
that the functor $\Ind:\u_\ell\mod\to \U_\ell\mod$ is exact. Therefore, 
it is sufficient to show that the map $\Res(\Ind(M))\to M$ is surjective when
when $M$ is an irreducible $\u_\ell$-module.

However, it is known (cf. \cite{Lu2} Proposition 5.11) that every irreducible
$\u_\ell$-module is of the form $\Res(N)$ for an irreducible $N\in \U_\ell\mod$.

Now, for $N$ as above the co-action map defines a map
$N\to \Ind(\Res(N))$, and the composition
$$\Res(N)\to \Res(\Ind(\Res(N)))\to \Res(N)$$
is the identity map.

Hence, 
$\Res(\Ind(M))\simeq \Res(\Ind(\Res(N)))\to \Res(N)\simeq M$
is a surjection.

\end{proof}

Condition (i) of \secref{coalset}
follows immediately from the fact that $\u_\ell$ acts trivially on any module of the
form $\Fr^*(V)$ for $V\in\check G\mod$. To verify condition (ii) we will use \propref{factorization}:

Indeed, we know that if $M$ is a $\Ag$-comodule,
then the co-action map restricted to $M^{\ag}$ factors as
$$M^{\ag}\to \Og\otimes M^{\ag}\hookrightarrow \Ag\otimes M.$$

By taking $M=\Ag$ and evaluating $(\on{id}\otimes \epsilon)\circ \Delta$
on $a\in \Ag^{\ag}$, we obtain that 
$$a=(\on{id}\otimes \epsilon)\circ \Delta(a)\in \Og.$$

\medskip

Let us now verify (iii). This follows from the next proposition:

\begin{prop}  \label{the ideal}
Let $(\CA,\CO,a)$ be satisfying properties (i), (ii). Suppose that
the adjunction map $\Res(\Ind(M))\to M$ is surjective for any 
$\a$-comodule $M$. Then $\m\cdot \CA=\on{Ker}(\CA\to \a)$.
\end{prop}

\begin{proof}

Set $I'=\m\cdot \CA$, $I=\on{Ker}(\CA\to \a)$. Set $\a':=\CA/I'$. Then $\a'$
is a co-algebra and a right $\CA$-module, and we have a sequence of epimorphisms
$\CA\to \a'\to \a$
respecting both structures.

We must show that the inclusion $I'\subset I$ is an equality.
For this, it is enough to show that $I$ is $\a'$-stable, i.e.
that the composition 
$$I\to \CA\overset{\Delta}\to \CA\otimes\CA\to \a'\otimes \CA$$
maps to $\a'\otimes I$. Indeed, by applying 
$(\on{id}\otimes\epsilon)\circ\Delta$ to $a\in I$ we then obtain that
$a=(\on{id}\otimes\epsilon)\circ\Delta(a)$ projects to the $0$ element in
$\a'$, i.e. belongs to $I'$.

\medskip

Using the fact that $\Res(\Ind(M))\to M$ i surjective for any $\a$-comodule,
we can find an $\CA$-comodule $B$ with a surjection $N\twoheadrightarrow I$.

However, from condition (ii), we obtain that in
$$\CA\hookrightarrow \CA^{\a'}\hookrightarrow \CA^{\a}$$
the composition is an isomorphism. Hence $\CA^{\a'}=\CA^{\a}$,
which implies that $N^{\a'}=N^{\a}$ for any $N\in\CA\comod$.
In particular, if $N_1$ and $N_2$ are two $\CA$-comodules,
any map $N_1\to N_2$ respecting the $\a$-coaction, respects
also the $\a'$-coaction.

Applying this to the composition
$N\twoheadrightarrow I\hookrightarrow \CA$, we obtain that 
$I=\on{Im}(N)\subset\CA$ is an $\a'$-subcomodule.

\end{proof}

Finally, as was mentioned above, the functor $\Ind$ is exact, 
hence $(\Ag,\Og,\ag)$ verifies condition (iv b).

For completeness, we will show that in fact $(\Ag,\Og,\ag)$ verifies also
condition (iv a). More generally, we will prove the following proposition:

\begin{prop}  \label{flatness}
let $\CO\to \CA$ be an embedding of Hopf algebras, with
$\CO$ being commutative. Then $\CA$ is faithfully-flat
as an $\CO$-module.
\end{prop}

\ssec{Proof of \propref{flatness}}

Let us denote $\on{Spec}(\CO)$ by $\Gamma$, and view $\CA$
as a quasi-coherent sheaf on $\Gamma$.

We have the following lemma:

\begin{lem}  \label{allequal}
For every $\gamma\in\Gamma$, the pull-back $\gamma^*(\CA)$ of $\CA$
under the translation map $\gamma'\to \gamma\cdot \gamma'$ is
(non-canonically) isomorphic to $\CA$ as a quasi-coherent sheaf.
\end{lem}

\begin{proof}
Let $\gamma\in\Gamma$ be a point over which the embedding
$\CO\to \CA$ induces an injection on fibers $\CO_\gamma\to \CA_\gamma$.
We will call such $\gamma$'s ``good''. First, we claim that for a ``good''
$\gamma$ we do obtain an isomorphism
$$\gamma^*(\CA)\simeq \CA.$$
Indeed, let $\xi_\gamma$ be any linear functional
$(\CA)_\gamma\to\CC$ which extends the evaluation map $\CO\to\CC$
corresponding to $\gamma$. Consider the map
$$\CA\overset{\Delta}\to \CA\ten\CA\overset{\xi_\gamma\ten\on{id}}
\longrightarrow \CA.$$
It is easy to see that this map defines the sought-for isomorphism
$\gamma^*(\CA)\simeq \CA$.

\medskip

Now let us show that all $\gamma\in \Gamma$ are ``good''. Suppose not.
Since $\CO\to \CA$ is an embedding , there exists a collection
$\underset{k}\cup Y_k$ of proper sub-schemes of $\Gamma$
defined over $\CC$, such that all points in $\Gamma\setminus \underset{k}\cup Y_k$
are ``good''. By what we proved above, the translation by a ``good''
$\gamma$ maps the collection $\underset{k}\cup Y_k$ to itself.

Let us make a field
extension $\CC\mapsto \CC(\Gamma)$. Over this field, $\Gamma$ has the
canonical generic point, which is clearly ``good''. However, this generic
point cannot map a proper sub-scheme defined over $\CC$ to another proper
sub-scheme defined over $\CC$, which is a contradiction.
\end{proof}

This lemma implies \propref{flatness}:

To prove that $\CA$ is flat over $\CO$, 
we must show that $\on{Tor}_{\CO}^1(\CA,\CC_\gamma)=0$, for every
$\gamma\in \Gamma$. (Here $\CC_\gamma$ denotes the sky-scraper sheaf at $\gamma$.)
As in the above argument, $\on{Tor}_{\CO}^1(\CA,\CC_\gamma)=0$ for all
$\gamma$'s lying outside $\underset{k}\cup Y_k$.

However, by \lemref{allequal},
all points of $\Gamma$ are ``the same'' with respect to $\CA$. Hence,
$\on{Tor}_{\CO}^1(\CA,\CC_\gamma)=0$ everywhere.

\medskip

To complete the proof of the proposition,
we must show that the fiber of $\CA$ at any $\gamma\in \Gamma$ is
non-zero. But this has been established in the course of the proof
of \lemref{allequal}.

\section{Further properties of the equivalence of categories}

\ssec{Definition by the universal property}

In this section we will make several additional remarks about the
equivalence of categories established in \thmref{main}. When our
discussion applies to any triple $(\CO,\CA,\a)$, we will work in this
more general context.
Let us denote by $F^*$ the natural functor from $\CO\comod$ to $\CA\comod$.

\medskip

By condition (i) of \secref{coalset}, we have an isomorphism of functors
$$\CO\comod\times \CA\comod\to \a\comod:\, 
\alpha_V^{can}:\Res(F^*(V)\otimes N)\simeq \underline{V}\otimes \Res(N).$$

\medskip

Let $\C$ be an abelian $\CC$-linear category and let
$\R:\CA\comod\to \C$ be a $\CC$-linear functor with the
property that for each $V\in \CO\comod$ and $N\in\CA\comod$ there
is a natural transformation
$$\alpha_V^{\C}: \R(F^*(V)\ten N) \mapsto
\underline{V}\ten \R(F^*(N)),$$
which satisfies the three properties of \secref{reform}.

\begin{prop} \label{universe}
There exists a functor
$\rr:\a\comod\to \C$ and an isomorphism of functors
$\R\simeq \rr\circ \Res$, such that $\alpha_V^{\C}=\rr(\alpha_V^{can})$.
\end{prop}

The meaning of this proposition is, of course, that the forgetful functor
$\Res:\CA\comod\to\a\comod$ is universal with respect to the property that
it transforms $F^*(V)\ten N$ to $\underline{V}\ten N$.

\begin{proof}

Using \thmref{gener}, we will think of $\a\comod$ in terms of $\Cat$
and we will construct a functor $\rr:\Cat\to\C$.

Let $\CO\ten M\overset{\on{act}}\to M$ be the action map.
By assumption, we obtain the map $$\underline{\CO}\ten \R(M)\simeq
\R(\CO\ten M)\to \R(M).$$
The axioms on the $\alpha_V^{\C}$'s imply that $\CO$ acts on $\R(M)$ as an
associative algebra. We set $\rr(M):=\CC\underset{\CO}\otimes \R(M)$.

Let us show now that $\R$ is canonically isomorphic to $\rr\circ \Res$.
Under the equivalence of
\thmref{main}, the functor $\Res$ goes over to $N\mapsto \CO\otimes N$.
Therefore,
$$\rr\circ \Res(N)\simeq 
\CC\underset{\CO}\otimes \R(\CO\otimes N)\simeq \R(N).$$
\end{proof}

\ssec{Reconstruction of $\CA\comod$ from $\a\comod$}

For $(\CO,\CA,\a)$ with $\CO$ being commutative,
let us recall from \corref{dual group action}
that the group $\Gamma=\Spec(\CO)$ acts on
the category $\a\comod$ by endo-functors.

Thus, it makes sense to talk about $\Gamma$-equivariant objects of
$\a\comod$.

\begin{prop} \label{backequiv}
The category of $\Gamma$-equivariant objects of $\a\comod$ is
naturally equivalent to $\CA\comod$.
\end{prop}

\begin{proof}
Let $M$ be a $H$-equivariant object in $\Cat$.
By definition, the underlying $\CA$-comodule has an additional
commuting structure of a $\Gamma$-equivariant quasi-coherent sheaf
on $\Gamma$. By taking its fiber at the point $1\in \Gamma$, we obtain an
$\CA$-comodule.

Thus, we have constructed a functor
$$\Cat_G\to\CA\comod,$$
and it is easy to see that it is an equivalence.
\end{proof}

Thus, given $\CA$, the category of $\a$-comodules is a
``de-equivariantization'' of $\CA\comod$.

\ssec{Other versions of quantum groups}

Let us discuss briefly the generalization of \thmref{main}
in the context of $\u_{\ell,sc}$ and $\tu_\ell$.

Consider the triple $\CA=\Ag$, $\CO=\CO_{\check G_{sc}}$ and
$\a=\ag{}_{sc}$. In a way completely analogous to what we did
in the previous section, one shows that these co-algebras satisfy
conditions (i)-(iii) of \secref{coalset}.

Let $\Catsc$ denote the corresponding category $\Cat$. We 
have the following version of \thmref{main}:

\begin{thm}   \label{main sc}
The categories $\Catsc$ and $\ag{}_{sc}\comod$ are naturally 
equivalent.
\end{thm}

\medskip

Now let us consider the case of $\tu_\ell$.
In what follows, for a $\check G$-module $V$, we will regard
$\Res^{\check G}_{\check T}(V)$ as a $Y$-graded vector space.

We introduce the category $\Catgd$ as follows: its objects are
$Y$-graded $\Ag$-comodules $M=\underset{\nu\in Y}\oplus M^\nu$,
each endowed with a collection of grading-preserving maps
$\alpha_V,\,V\in \check G\mod$
$$\Fr^*(V)\ten M\simeq \Res^{\check G}_{\check T}(V)\ten M,$$
(as in \secref{reform})
where the $Y$-grading on the LHS comes from the grading on $M$ and on
the RHS the grading is diagonal.
Maps in this category are grading preserving $\U_\ell$-module maps,
which intertwine the corresponding $\alpha_V$'s.

\begin{thm} \label{maingraded}
The category $\Catgd$ is equivalent to $\tu_\ell\mod$.
\end{thm}

\begin{proof}

First, let us observe that if we put $\CA=\agd$, $\CO={\mathcal O}_{\check T}$,
$\a=\ag$, the corresponding triple would satisfy conditions (i)-(iii)
of \secref{coalset}. Hence, the general \thmref{gener} is applicable
as well as \propref{backequiv}.

Therefore, the category $\Catgd$ is equivalent to the category
of $\check T$-equivariant objects in $\Catg$. However, the latter
is by definition the same as $\Catgd$.

\end{proof}

\medskip 

Finally, let us characterize the category $\Catgd$ by a universal
property.

For an $\check\lambda$ in $Y=X^*$ (or even in $X^*_{sc}$) let us denote
by $\CC^{\check\lambda}$ the corresponding $1$-dimensional module over $\tu_\ell$,
and by $P_{\check\lambda}:\tu_\ell\mod\mapsto \tu_\ell\mod$ the translation functor
$M\mapsto \CC^{\check\lambda}\otimes M$.

Let now $\C$ be an abelian $\CC$-linear category and let $
P^{\C}_{\check\lambda}:\C\to \C$
be an action of $Y$ on $\C$ by endo-functors.
Let $\R:\Ag\comod\to \C$ be a $\CC$-linear functor with the property
that for each $V\in \Og\comod$ there is a natural transformation
$$\alpha_V^{\C}: \R(\Fr^*(V)\ten M) \mapsto \underset{\nu}\oplus \,
\underline{V}(\check\nu)\otimes P_{\check\nu}^{\C}(\R(M)),$$
which satisfies the three properties of \secref{reform}.
(In the above formula, for a $\check G$-module $V$ and $\check\nu\in Y$,
$\underline{V}(\check\nu)$ denotes the corresponding weight subspace.)

\begin{prop}
There exists a functor $\rr:\agd\comod\to \C$ and an isomorphism of functors
$$\R\simeq \rr\circ \Res.$$
Moreover, the functor $\rr$ commutes with the translation functors
in the obvious sense.
\end{prop}

We omit the proof, since it is completely analogous to the proof
of \propref{universe}.

\section{The regular block}  \label{regular block}

\ssec{Blocks in the categories ${\mathcal A}\operatorname{-comod}$
and ${\mathfrak a}\operatorname{-comod}$.} \label{reg}

Recall that any Artinian abelian category $\C$ is a direct sum of
its indecomposable abelian sub-categories called {\em blocks}
or {\it linkage classes} of $\C$. Obviously, a block of a category 
is completely described by the set of irreducible objects contained in it.
We will denote the set of blocks of $\C$ by $\Bl(\C)$.

Note that the categories of finite dimensional  $\CA$- and
$\a$-comodules (denoted below by $\CA\comod^{f}$ and $\a\comod^{f}$,
respectively) are Artinian, therefore, they admit decompositions into blocks.
We will use the notation $\Bl(\CA)$ and $\Bl(\a)$ for the sets of blocks 
of $\CA\comod^f$ and $\a\comod^f$, respectively.

Evidently, we have $$\CA\comod= \on{ind.comp.}
\left(\CA\comod^{f}\right)
\text{ and } \a\comod=\on{ind.comp.}\left(\a\comod^{f}\right).$$
For $\alpha\in \Bl(\CA)$ (resp., $\alpha'\in \Bl(\a)$)
let us denote by $\CA\comod_\alpha$ (resp., $\a\comod_{\alpha'}$)
the ind-completion of the corresponding block of $\CA\comod^f$
(resp., $\a\comod^{f}$).

\smallskip

We will call the block of
$\CA\comod$ (resp., $\a\comod$) which contains the trivial representation
$\CC$ {\em the regular block} and will denote it by $\CA\comod_0$
(resp., $\a\comod_0$).

\medskip

Assume that the category $\CA\comod$ has the following
additional property with respect to $\CO\comod$:
 
\smallskip

\noindent ($*$) {\it  For any $\alpha\in \Bl(\CA)$ and 
$V\in\CO\comod$, the functor
$F^*(V)\ten\cdot:\ \CA\comod\map\CA\comod$
maps $\CA\comod_\alpha$ to itself.}

\bigskip

Let us compare the block decompositions of $\CA\comod$ and  $\a\comod$.

\begin{prop} \label{onetoone}
There is a one-to-one correspondence between the sets $\Bl(\CA)$
and $\Bl(\a)$ determined by the following properties:

\smallskip
\noindent
{\em (a)} $N\in \CA\comod_\alpha$ if and only if $\Res(N)\in
\a\comod_\alpha$.

\smallskip
\noindent {\em (b)}
$M\in\a\comod_\alpha$ if and only if $\Ind(M)\in\CA\comod_\alpha$.

\end{prop}

\begin{proof}
First, observe that $\Ind\circ\Res:\CA\comod\to\CA\comod$
preserves each $\CA\comod_\alpha$, by assumption, since $\Ind\circ\Res(N)
\simeq F^*(\CO)\ten N$.

Secondly, let us show that $\Res\circ\Ind$ maps
each $\a\comod_\alpha$ to itself. Indeed, let $M\in \a\comod_\alpha$
and let $N'$ be an $\a$-stable direct summand of $\Ind(M)$,
which belongs to some $\a\comod_\beta$. Then $N'$ is preserved by the
$\CO$-action, and thus defines a sub-object of $\underline\Ind(M)\in\Cat$.
But then $\underline\Res(N')\in\a\comod_\beta$ is a non-zero direct summand 
of $M$, which means that $\beta=\alpha$.

\medskip

Let $N$ be an object of $\CA\comod$. Let
$\Res(N)=\Res(N)'\oplus \Res(N)''$ be a block decomposition in $\a\comod$.
Let us show that $\Res(N)'$ and $\Res(N)''$ are in fact $\CA$-sub-comodules.
Without restricting the generality, we can assume that $N$ is a sub-comodule
of $\Ind(M)$ for some $M\in \a\comod$. However, as we have just seen,
the block decomposition of $\Res\circ\Ind(M)$
coincides with the block decomposition of $M$.

\medskip

Therefore, the block decomposition of $\a\comod$ is ``coarser'' than
that of $\CA\comod$.

However, by our assumption on $\CA$, its block decomposition
is ``coarser'' than the block decomposition of
$\Cat$. This implies the assertion of the proposition
in view of \thmref{reform}.
\end{proof}

\ssec{The category $\Cat^0$} \label{newprop}

We define the category $\Cat^0$ as the preimage of $\CA\comod_0$ under
the tautological forgetful functor $\Cat\to\CA\comod$. This definition
makes sense due to condition $(*)$ above.
In the course of the proof of \propref{onetoone} we have established the
following assertion:

\begin{cor} \label{absregular}
Under the equivalence of categories $\Cat\simeq\a\comod$,
the sub-category $\Cat^0$ goes over to the regular block
$\a\comod_0$.
\end{cor}

\ssec{The case of $\u_\ell$}
For a regular dominant $\lambda\in X$, let $\W(\lambda)\in \U_\ell\mod$
denote the corresponding {\it Weyl module}. It is well-known that
$\W(\lambda)$ has a unique simple quotient, denoted $\L(\lambda)$ and that
each simple object in $\U_\ell\mod$ is isomorphic to $\L(\lambda)$
for some $\lambda$.

The following facts about the block decomposition of 
the category $\U_\ell\mod$ were established in \cite{APW}:

Let $W_{aff}$ be the affine Weyl group $Y\ltimes W$.
It acts on the lattice $X$ as follows: the translations by $Y$ act via
the homomorphism $\phi:Y\to X$ and the action of the finite Weyl group
$W$ is centered at $-\rho:=-\underset{i\in I}\Sigma\, \omega_i$, 
where $\omega_i$'s are the fundamental weights.

\begin{thm} \label{APW}
Two simple
modules $\L(\lambda_1)$ and $\L(\lambda_2)$ are in  the same  block
if and only if $\lambda_1$ and $\lambda_2$ belong to the same
$W_{aff}$-orbit.
\end{thm}

Moreover, we have the following statement (cf. \cite{Lu4}, Theorem 7.4
and Proposition 7.5):

\begin{prop}  \label{irreducibles}
Let $\lambda=\lambda_1+\lambda_2$ be the unique decomposition,
with $\lambda_2=\phi_{sc}(\check\mu)$, where $\check\mu\in X^*_{sc}$
is a dominant integral weight of the group $\check G_{sc}$
and $\lambda_1$ is such that
$0\leq\langle \lambda_1,\check \alpha_i\rangle < \ell_i$ for all $i\in I$. 
Then:
\begin{itemize}

\item[(i)]
$\L(\lambda_2)\simeq \on{Fr}_{sc}^*(V^{\check\mu})$,
where $V^{\check\mu}$ is the corresponding irreducible representation of
$\check G_{sc}$.

\item[(ii)]
The restriction of $\L(\lambda_1)$ to $\u_\ell$ remains irreducible.

\item[(iii)]
$\L(\lambda)\simeq \L(\lambda_1)\otimes \L(\lambda_2)$.
\end{itemize}
\end{prop}

This proposition combined with \thmref{APW} implies that the 
category $\U_\ell\mod$ satisfies condition $(*)$.

\medskip

Let $\Catg^0$ denote the corresponding sub-category of $\Catg$.
By applying \propref{onetoone} and \corref{absregular},
we obtain the following theorem:

\begin{thm} \label{mainregular}
We have a bijection between the sets $\Bl(\U_\ell\mod)\simeq
\Bl(\u_\ell\mod)$ and an equivalence of categories:
$$\u_\ell\mod_0\simeq \Catg^0.$$
\end{thm}

Recall (cf. \cite{Lu2} or \cite{APW1}) that to every element
$\lambda\in X$ we attached an irreducible
object of $\u_\ell\mod$, denoted $L(\lambda)$, which depends only
on the image of $\lambda$ in the quotient $X/\phi(Y)$, and
\propref{irreducibles}(ii) implies that if $\lambda$ satisfies
$\langle \lambda,\check \alpha_i\rangle < \ell_i$, then
$L(\lambda)\simeq \Res({\mathbf L}(\lambda))$.

The following corollary repeats in fact Sect. 2.9 of \cite{APW1}:

\begin{cor}
For two elements $\lambda_1$ and $\lambda_2$ of $X$, the modules
$L(\lambda_1)$ and $L(\lambda_2)$
belong to the same block of $\u_\ell\mod$ if and only if $\lambda_1$
and $\lambda_2\in X$ are $W_{aff}$-conjugate. 
\end{cor}

\ssec{The graded case}

For completeness, let us analyze the block decomposition
of the category $\tu_\ell\mod$ in light of \thmref{onetoone}.
However, all that we are going to obtain is already contained
in \cite{APW1}.

Recall that for $\lambda\in X$, $\TL(\lambda)$ denotes
the corresponding irreducible object of $\tu_\ell\mod$, 
and $L(\lambda)=\Res^{\tu_\ell}_{\u_\ell}(\TL(\lambda))$.

The following is known, due to \cite{APW1}:

\begin{prop}
The translation functors $P_{\check\lambda}$, 
$\check\lambda\in Y$ preserve the block
decomposition of $\tu_\ell\mod$.
\end{prop}

Using this proposition, we can apply \propref{onetoone} to the category
$\tu_\ell\mod$ and the group $\check T$. Thus, we obtain
the following result of \cite{APW1}:

\begin{cor}
There is a natural bijection $\Bl(\tu_\ell\mod)\simeq \Bl(\u_\ell\mod)$.
The modules $\TL(\lambda_1)$ and $\TL(\lambda_2)$ belong to the block in
$\tu_\ell\mod$ if and only if $\lambda_1$ and $\lambda_2$ are
$W_{aff}$-conjugate. 
\end{cor}

Let $\Catgd^0$ denote the preimage of $\U_\ell\mod_0$ under the obvious
forgetful functor. From \propref{onetoone} and \thmref{maingraded}, we
obtain the following theorem (cf. \cite{AJS} for the first assertion):

\begin{thm} \label{mainregulargraded}
There is an isomorphism of sets $\Bl(\U_\ell\mod)\simeq \Bl(\tu_\ell\mod)$
and an equivalence of categories:
$$\tu_\ell\mod_0\simeq \Catgd^0.$$
\end{thm}

\ssec{The case of $\u_{\ell,sc}$.}

Observe that if we consider the triple $\CA=\Ag$, $\CO=\CO_{\check G_{sc}}$,
$\a=\ag{}_{sc}$, then condition (*) above will not be satisfied.
Instead, we have the following assertion:

\begin{prop}
The natural restriction functor $\u_\ell\mod\to \u_{\ell,sc}\mod$
induces an equivalence $\u_\ell\mod_0\to \u_{\ell,sc}\mod_0$.
\end{prop}

\begin{proof}

For $\lambda\in X$, let us denote by $L(\lambda)_{sc}$ the restriction of
$L(\lambda)$ to $\u_{\ell,sc}$. By construction, it depends only
on the class of $\lambda$ in $X/\phi(X^*_{sc})$. 

Let us consider the forgetful functor 
$\Res^{\u_\ell}_{\u_{\ell,sc}}:\u_\ell\mod\to\u_{\ell,sc}\mod$.
Note that in terms of $\Catg$ and $\Catsc$, it acts as follows:
$$M\in \Catg\mapsto \CO_{\check G_{sc}}\underset{\Og}\otimes M\in \Catsc.$$ 
This functor has a right adjoint, which we will denote by
$\Ind^{\u_\ell}_{\u_{\ell,sc}}$.
On the level of $\Catsc$, $\Ind^{\u_\ell}_{\u_{\ell,sc}}$
is the natural forgetful functor.

Note that for $M\in \u_\ell\mod$ we have:
$$\Ind^{\u_\ell}_{\u_{\ell,sc}}\circ \Res^{\u_\ell}_{\u_{\ell,sc}}(M)\simeq
\underset{\check\lambda\in X^*_{sc}/X^*}\oplus\, M\otimes \CC^{\check\lambda}.$$
(Recall that $\CC^{\check\lambda}$ is a $1$-dimensional $\tu_\ell$-module,
and, hence, we are allowed to tensor any $\u_\ell$-module by it on the right.)
In particular, $\Ind^{\u_\ell}_{\u_{\ell,sc}}(L(\mu)_{sc})=
\underset{\check\lambda\in X^*_{sc}/X^*}\oplus\, L(\mu+\phi_{sc}(\check\lambda))$.

Therefore, two irreducible objects $L(\mu^1)_{sc}$ and 
$L(\mu^2)_{sc}$ of $\u_{\ell,sc}\mod$ belong to the same block
if and only if there exists $\check\lambda\in X^*_{sc}$, such that 
$L(\mu^1+\phi_{sc}(\check\lambda))$ and $L(\mu^2)$ belong to the same block of $\u_\ell\mod$,
i.e. $\mu^1$ and $\mu^2$ belong to the same orbit of the extended
affine Weyl group $W^{ext}_{aff}\simeq X^*_{sc} \ltimes W$.

\medskip

Thus, we obtain that the functor $\Res^{\u_\ell}_{\u_{\ell,sc}}$ maps
$\u_\ell\mod_0$ to $\u_{\ell,sc}\mod_0$. We claim now that this functor 
has a left quasi-inverse. Namely, it is given by 
$$N\mapsto pr_0(\Ind^{\u_\ell}_{\u_{\ell,sc}}(N)),$$
where $pr_0$ denotes the functor of projection onto the regular block
in $\u_\ell\mod$. Indeed, for $M\in \u_\ell\mod_0$ we have:
$$pr_0(\Ind^{\u_\ell}_{\u_{\ell,sc}}\circ \Res^{\u_\ell}_{\u_{\ell,sc}}(M))\simeq M,$$
because for $\check\lambda\in X^*_{sc}$, the object 
$pr_0(M\otimes\CC^{\check\lambda})$ 
is non-zero only if $\check\lambda\in X^*$, which follows from the description 
of blocks of $\u_\ell\mod$ in terms of $W_{aff}$.

To finish the proof of the proposition it remains to show that if
$N$ is a non-zero object in $\u_{\ell,sc}\mod_0$, 
then $pr_0(\Ind^{\u_\ell}_{\u_{\ell,sc}}(N))$ is non-zero either. For that,
it is enough to suppose that $N$ is irreducible, i.e. of the form
$L(\lambda)_{sc}$, and our assertion follows from the explicit description
of $\u_{\ell,sc}\mod_0$ given above.

\end{proof}

\section{Geometric interpretation} \ssec{Affine flag variety}

Our goal in this section is to give a geometric description of
the category $\u_\ell\mod_0$.
Namely, we will show that it can be described as the category of
certain perverse sheaves on the affine flag variety corresponding to the
group $G$, which have the Hecke eigen-property.
It is via this description that one can link $\u_\ell\mod_0$ to
certain categories which appear in the geometric Langlands correspondence
and to other interesting categories arising in representation
theory.

\smallskip

First, we will briefly recall several definitions concerning the affine
Grassmannian and the affine flag variety. We refer the reader to 
\cite{MV}, \cite{BD} or \cite{Ga} for a more detailed discussion.

\medskip

Consider the ring $\CC[[t]]$ of Taylor series and the field $\CC((t))$
of Laurent series. The loop group $G((t))$ (resp., the group of positive
loops $G[[t]]$) has a structure of an indgroup-scheme (resp., of a
group-scheme). The quotient $G((t))/G[[t]]$ is an ind-scheme 
of ind-finite type, called {\it the
affine Grassmannian} of the group $G$, denoted $\Gr$.

\smallskip

By definition, the group-scheme $G[[t]]$ acts (on the left) on $\Gr$.
The orbits of this action are finite-dimensional quasi-projective 
varieties and
they are in a natural bijection with the dominant elements $Y^+\subset Y$.
For $\check\lambda\in Y^+$, we will denote by $\overline{\Gr}^{\check\lambda}$
the closure of the corresponding orbit. Thus, it makes sense to talk
about the category of $G[[t]]$-equivariant perverse sheaves on
$\Gr$. By definition, every such perverse sheaf is
supported on $\overline{\Gr}^{\check\lambda}$ for ${\check\lambda}$ sufficiently large.
We will denote this category by $\PG$. This is an abelian category and
it possesses an additional structure of the {\it convolution product}
$\PG\star\PG\to\PG$, which makes $\PG$ into a tensor category. 

We have the following fundamental theorem (\cite{Gi}, \cite{MV}):

\begin{thm} \label{Langlands}
There is an equivalence of tensor categories $\check G\mod\simeq \PG$. 
Under this equivalence,
the intersection cohomology sheaf $\on{IC}_{\overline{\Gr}^{\check\lambda}}$ goes
over to the highest weight module $V^{\check\lambda}$.
\end{thm}

Now we pass to the definition of the affine flag variety.

\medskip

Let $\Iw\subset G[[t]]$ be the Iwahori subgroup. By definition, $\Iw$ is
the preimage of the Borel subgroup $B\subset G$ under the natural evaluation
map $G[[t]]\to G$. The quotient $G((t))/\Iw$
is also an ind-scheme of ind-finite type,
called {\it the affine flag variety} of $G$,
denoted $\Fl$. By definition, we have a projection
$\Fl\twoheadrightarrow\Gr$,
whose fibers are (non-canonically) isomorphic to the usual flag manifold
$G/B$.

Let $N_{\Iw}$ be the unipotent radical of $\Iw$. By definition, $N_{\Iw}$ is
the preimage of $N\subset B$ under
$\Iw\to B$. Since $N_{\Iw}$ is normal in $\Iw$ and $\Iw/N_{\Iw}\simeq T$,
the quotient $\Flt:=G((t))/N_{\Iw}$ is a principal $T$-bundle over
$\Fl:=G((t))/\Iw$.
We  will call $\Flt$ the enhanced affine flag variety.
The group-scheme $G[[t]]$ acts on the left on both $\Fl$ and $\Flt$.

We define the category $\PF$ to be the abelian category of
$G[[t]]$-equivariant perverse sheaves on $\Fl$. We define the category
$\PFt$ to be the sub-category of the category of $G[[t]]$-equivariant
perverse sheaves on $\Flt$, which consists of $T$-monodromic objects.
\footnote{Recall that a perverse sheaf is called $T$-monodromic
when it has a filtration, whose sub-quotients are $T$-equivariant perverse sheaves.}
Note that the pull-back functor identifies  $\PF$ with the sub-category
of $\PFt$ consisting of $T$-equivariant objects.

\smallskip

It is known (cf. \cite{Ga} for details) that the 
convolution tensor structure on the
category $\PG$ extends to an action of $\PG$ on $\PF$.

Similarly, one can define an action $$\PG\star\PFt\to\PFt.$$

\ssec{Hecke eigen-sheaves}
Let $\PFti$ denote the ind-completion of the category $\PFt$.
We define the category $\Af$ as follows:
its objects are pairs 
$$(\F\in\PFti,\{\alpha_V,\,\,\forall\,V\in\check G\mod\}),$$
where each $\alpha_V$ is a map
$$\F_V\star \F\to \underline{V}\ten\F,$$ where $\F_V\in \PG$ is
the perverse sheaf corresponding to $V\in\check G\mod$ via
the equivalence of categories of \thmref{Langlands}.
The maps $\alpha_V$ must satisfy the three conditions of \secref{reform}.
Morphisms in $\Af$ between $(\F,\alpha_V)$ and $(\F',\alpha'_V)$ are
maps $\F\to \F'$, which intertwine between $\alpha_V$ and $\alpha'_V$.
As in \propref{isom}, one shows that the maps $\alpha_V$ as above are automatically isomorphisms.

\medskip

The rest of this section (and of the paper) is devoted to the proof of
the following theorem.

\begin{thm} \label{geominter}
For $\ell$ sufficiently large, there is an equivalence of categories between
$\Af$ and $\u_\ell\mod_0$.
\end{thm}

Unfortunately, the proof relies on two results, whose proofs
are unavailable in the published literature. Therefore, the reader may
regard \thmref{geominter} as a conjecture, which can be deduced from
\thmref{Finkelberg} and \thmref{convolution} stated below.

\ssec{Twisted D-modules on $\Flt$} \label{KLKT}

The first step in the passage $\Af\map\u_\ell\mod_0$ is
the functor from perverse sheaves on $\Flt$ to modules over the Kac-Moody
algebra due to \cite{KT}. 

Recall that to an invariant symmetric form $c:\g\otimes\g\to \CC$,
which is {\it integral} (i.e. induces an integral-valued form on 
the cocharacter lattice $Y$), we can associate a line bundle
${\mathcal L}_c$ on $\Gr$ (cf. \cite{KT}). By pulling
it back to $\Fl$ and $\Flt$ we obtain the corresponding
line bundles on the latter ind-schemes.

Thus, we can consider the category of
${\mathcal L}_c$-twisted right D-modules on $\Flt$, cf. \cite{KT},
\cite{BD}. As ${\mathcal L}_c$ is $G[[t]]$-equivariant, it makes sense to
consider the category $\DFlc$ of $G[[t]]$- equivariant $T$-monodromic
${\mathcal L}_c$-twisted right D-modules on $\Flt$.

\smallskip

Let us now consider the Kac-Moody algebra $\gh$ corresponding to $\g$:
$$0\to \CC\to \gh\to \g((t))\to 0,$$
defined with respect to the pairing $\g\otimes\g\to \CC$ given
by $c$. Let us denote by $\gh_c\mod$ the category of continuous representations
of $\gh$ on which $1\in \CC\subset \gh$ acts as identity.

Let us denote by $\gh_c^{G[[t]]}\mod$ the subcategory of $\gh_c\mod$
consisting of finite length representations,
on which the action of $\g[[t]]\subset \gh$ integrates to the action of the
group-scheme $G[[t]]$.
According to \cite{KL}, this is an Artinian category and we will
denote by $\KMc$ its regular block.

\medskip

From now on let us suppose that $c$ is such that $c_{crit}-c$ is 
positive definite on $Y$, where $c_{crit}$ corresponds to 
$-\frac{1}{2}\cdot$(the Killing form). When $\g$ is simple,
$c_{crit}$ is $-\check h$ times the {\it canonical} integral
from on $Y$, where $\check h$ is the dual Coxeter number.

The following theorem has been established in \cite{KT} and \cite{BD}:

\begin{thm}
The functor of global sections of a twisted D-module defines an exact and
faithful functor:
$$\DFlc\map\KMc.$$
\end{thm}

However, a stronger statement is true:
\footnote{This result is probably well-known to many experts.
The proof that we have in mind has been explained to us by M.~Finkelberg
and R.~Bezrukavnikov.}

\begin{thm} \label{Finkelberg}
The above functor
$\DFlc\map\KMc$ is in fact an equivalence of categories.
\end{thm}

The Riemann-Hilbert correspondence yields an equivalence
of categories between $\DFlc$ and $\PFt$, cf. \cite{KT}. Therefore, we
obtain the following corollary:

\begin{cor} \label{sheavesKM}
There is an equivalence of categories $\PFt\simeq\KMc$.
\end{cor}

\ssec{The Kazhdan-Lusztig equivalence of categories}

Now let $\ell$ be as in \secref{biggroup} and set
$c=c_{crit}-(\cdot,\cdot)_\ell$, where $(\cdot,\cdot)_\ell$ 
has been introduced in \secref{quantum Frob}. Again, when $\g$ is simple,
$c$ is $-\check h-\frac{\ell}{d}$ times the canonical form, where
$d=\on{max}(d_i)$.

The following theorem has been established in \cite{KL}:

\begin{thm}
When $\ell$ is sufficiently large, there is an equivalence of categories
$\U_\ell\mod\simeq \gh_c^{G[[t]]}\mod$.
\end{thm}

By combining this theorem with \corref{sheavesKM}, we obtain the
following corollary:

\begin{cor}
\label{sheavesquant}
There is an equivalence of categories: $\PFt\simeq
\U_\ell\mod_0$.
\end{cor}

To prove \thmref{geominter}, we will need the following
property of the equivalence stated in \corref{sheavesquant}:

\begin{thm}
\label{convolution}
Under the above equivalence of categories
$\PFt\simeq \U_\ell\mod_0$
the functors
$\check G\mod\times \PFt\to\PFt$ given by
$$V,\F\mapsto \F_V\star \F \text{ and } V,M\mapsto \Fr^*(V)\ten M$$
are naturally isomorphic.
\end{thm}

This result has not been stated explicitly in \cite{KL}.
We will supply the proof in a later publication.

\smallskip

Now, by passing to the ind-completions of the categories $\PFt$ and
$\U_\ell\mod_0$, we obtain that \thmref{geominter} is a
consequence of \thmref{main} and \thmref{convolution}.

\end{document}